\documentclass[24pt]{amsart}
\usepackage{amsmath, amsthm, amscd, amsfonts}
\usepackage{graphicx}
\usepackage{graphics}

\newtheorem{thm}{Theorem}[section]

\newtheorem{prop}[thm]{Proposition}
\newtheorem{cor}[thm]{Corollary}
\theoremstyle{definition}

\newtheorem{note-prob}[thm]{Note and Problem}

\theoremstyle{remark}
\newtheorem{rem}[thm]{Remark}
\numberwithin{equation}{section}
\newfont{\kh}{msbm10}

\def\R{{\mathbb R}}
\def\S{\mathfrak{S}}
\def\H{\mathfrak{H}}
\def\F{\mathfrak{F}}
\def\T{\mathfrak{T}}
\def\P{\mathfrak{P}}

\begin{document}

\title{Identities by Generalized $L-$Summing Method}

\author{M. Hassani}
\address{Mehdi Hassani, \newline Department of Mathematics, Institute for Advanced Studies in Basic
Sciences, P.O. Box 45195-1159, Zanjan, Iran}
\email{mmhassany@member.ams.org}

\author{Z. Jafari}
\address{Zahra Jafari, \newline Department of Information Technology, Institute for Advanced Studies in Basic
Sciences, P.O. Box 45195-1159, Zanjan, Iran}
\email{z.jafari@iasbs.ac.ir, zahra7jfr@yahoo.com}


\begin{abstract}
In this paper, we introduce 3-dimensional $L-$summing method, which
is a rearrangement of the summation $\sum A_{abc}$ with $1\leq
a,b,c\leq n$. Applying this method on some special arrays, we obtain
some identities on the Riemann zeta function and digamma function.
Also, we give a Maple program for this method to obtain identities
with input various arrays and out put identities concerning some
elementary functions and hypergeometric functions. Finally, we
introduce a further generalization of $L-$summing method in higher
dimension spaces.
\end{abstract}

\maketitle

\section{Introduction and Motivation}
Consider the following $n\times n$ multiplication table
\begin{center}
\begin{figure}[ht]
 \vspace{0cm} \hspace{0cm}
 \includegraphics[height=4cm,width=8cm]{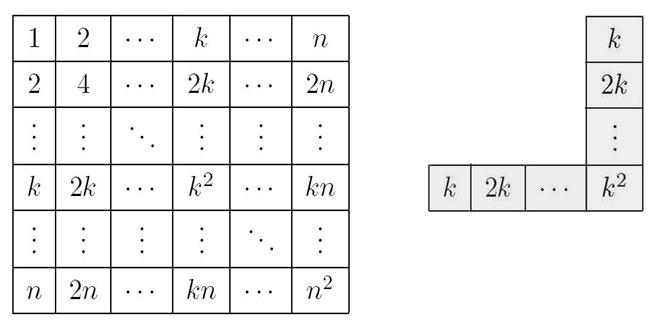}
 \caption{\Small{Multiplication table and $L-$summing element, $L_k$}}
\end{figure}
\end{center}
If we set $\Sigma(n)$ for the sum of all numbers in it, then by
summing line by line we have
$\Sigma(n)=\left(\frac{n(n+1)}{2}\right)^2$. On the other hand, we
can find $\Sigma(n)$ by using another method; letting $L_k$ be the
sum of numbers in the rotated $L$ in above table (right part of
\textsc{Figure} 1), we have
$$
L_k=k+2k+\cdots+k^2+\cdots+2k+k=2k(1+2+\cdots+k)-k^2=k^3.
$$
We call $L_k$, $L-$summing element. Thus we get
$\Sigma(n)=\sum_{k=1}^n L_k=\sum_{k=1}^n k^3$, and therefore
$\sum_{k=1}^n k^3=\left(\frac{n(n+1)}{2}\right)^2$. This is
2-dimensional \textit{L-}summing method (applied on the array
$A_{ab}=ab$), which briefly is
\begin{equation}\label{l-summing}
\sum(L-{\rm Summing~Elements})=\Sigma.
\end{equation}
More precisely, the $L-$summing method of elements of $n\times n$
array $A_{ab}$ with $1\leq a,b\leq n$, is the following
rearrangement
$$
\sum_{k=1}^n\left\{\sum_{a=1}^k A_{ak}+\sum_{b=1}^k
A_{kb}-A_{kk}\right\}=\sum_{1\leq a,b\leq n}A_{ab}.
$$
This method allows to obtain easily some classical algebraic
identities and also, with help of MAPLE, some new compact formulas
for sums related with the Riemann zeta function, the gamma function
and the digamma function \cite{jg, hassani}.\\ In this paper we
introduce a 3-dimensional version of $L-$summing method for $n\times
n\times n$ arrays and applying it on some special arrays we obtain
some identities concerning the Riemann zeta function and digamma
function. Then, we give a Maple program for this method and using it
we generate and then proof some new identities, concerning some
elementary functions and hypergeometric functions. Finally, we
introduce a further generalization of $L-$summing method in higher
dimension spaces and for latices related by a manifold.

\section{Formulation of the $L-$summing method in $\mathbb{R}^3$} Consider a three dimensional array $A_{abc}$ with
$1\leq a,b,c\leq n$ and $n$ is a positive integer. We should prepare
an explicit version of the general formulation (\ref{l-summing}) for
this array. The summation of all entries is $\Sigma(n)=\sum_{1\leq
a,b,c\leq n}A_{abc}$. The $L-$summing elements in this array have
the form pictured bellow
\begin{center}
\begin{figure}[ht]
 \vspace{0cm} \hspace{0cm}
 \includegraphics[height=4.5cm,width=6.5cm]{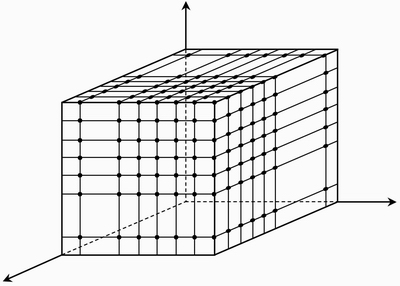}
 \caption{\Small{$L-$summing elements in $\mathbb{R}^3$}}
\end{figure}
\end{center}
So, we have $L_k=\Sigma_2-\Sigma_1+\Sigma_0$, with
$$
\Sigma_2=\sum_{b,c=1}^k A_{kbc}+\sum_{a,c=1}^k
A_{akc}+\sum_{a,b=1}^k A_{abk},\hspace{10mm}\Sigma_1=\sum_{a=1}^k
A_{akk}+\sum_{b=1}^k A_{kbk}+\sum_{c=1}^k
A_{kkc},\hspace{10mm}\Sigma_0=A_{kkk}.
$$
Therefore, $L-$summing method in $\mathbb{R}^3$ take the following
formulation
\begin{equation}\label{3d-l-summing}
\sum_{k=1}^n\{\Sigma_2-\Sigma_1+\Sigma_0\}=\Sigma(n).
\end{equation}
Note that $\Sigma_2$ is the sum of entries in three faces,
$\Sigma_1$ is the sum of entries in three intersected edges and
$\Sigma_0$ is
the end point of all faces and edges.\\
If the array $A_{abc}$ is symmetric, that is for each permutation
$\sigma\in S_3$ it satisfies
$A_{abc}=A_{\sigma(a)\sigma(b)\sigma(c)}$, then $L-$summing elements
in $\mathbb{R}^3$ takes the following easier form
\begin{equation}\label{symmetric-lk}
L_k=3\sum_{1\leq b,c\leq k} A_{kbc}-3\sum_{a=1}^k A_{akk}+A_{kkk}.
\end{equation}
In the next section we will apply 3-dimensional $L-$summing method
on two special symmetric arrays, related by the Riemann zeta
function and digamma function.

\section{Arrays related by the Riemann zeta function and digamma function}
\subsection{The Riemann zeta function} Suppose $s\in\mathbb{C}$ and let $A_{abc}=(abc)^{-s}$. It is
clear that
$$
\Sigma(n)=\sum_{1\leq a,b,c\leq
n}(abc)^{-s}=\left(\sum_{k=1}^n\frac{1}{k^s}\right)^3=\zeta_n^3(s),
$$
where $\zeta_n(s)=\sum_{k=1}^n k^{-s}$. Since this array is
symmetric, considering (\ref{symmetric-lk}), we have
$$
L_k=3\frac{\zeta_k^2(s)}{k^s}-3\frac{\zeta_k(s)}{k^{2s}}+\frac{1}{k^{3s}},
$$
and an easy simplifying, we can reform $\sum L_k=\Sigma(n)$ as
follows
\begin{equation}\label{l-summing-zeta-n}
\sum_{k=1}^n
\frac{\zeta_k^2(s)}{k^s}-\frac{\zeta_k(s)}{k^{2s}}=\frac{\zeta_n^3(s)-\zeta_n(3s)}{3}.
\end{equation}
Note that if $\Re(s)>1$, then
$\lim_{n\rightarrow\infty}\zeta_n(s)=\zeta(s)$, where
$\zeta(s)=\sum_{k=1}^{\infty}n^{-s}$ is the well-known Riemann zeta
function defined for complex values of $s$ with $\Re(s)>1$ and
admits a meromorphic continuation to whole complex plan
\cite{ivic-zeta}. So, for $\Re(s)>1$ we have
$$
\sum_{k=1}^{\infty}
\frac{\zeta_k^2(s)}{k^s}-\frac{\zeta_k(s)}{k^{2s}}=\frac{\zeta^3(s)-\zeta(3s)}{3},
$$
which also is true for other values of $s$ by meromorphic
continuation, except $s=1$ and $s=\frac{1}{3}$.

\subsection{Digamma function} Setting $s=1$ in (\ref{l-summing-zeta-n}) (or equivalently taking $A_{abc}=\frac{1}{abc}$)
and considering $\zeta_n(1)=H_n=\sum_{k=1}^n\frac{1}{k}$, we obtain
$$ \sum_{k=1}^n
\frac{H_k^2}{k}-\frac{H_k}{k^2}=\frac{H_n^3-\zeta_n(3)}{3}.
$$
One can state this identity in sense of digamma function
$\Psi(x)=\frac{d}{dx}\ln\Gamma(x)$, with
$\Gamma(x)=\int_0^{\infty}e^{-t}t^{x-1}dt$ is the well-known gamma
function. Considering logarithmic derivative of the formula
$\Gamma(n+1)=n\Gamma(n)$, we obtain
\begin{equation}\label{psi-n+1-n}
\Psi(n+1)=\frac{1}{n}+\Psi(n),
\end{equation}
and applying this relation, we yield that $\Psi(n+1)-\Psi(1)=H_n$.
Thus, we have
\begin{equation}\label{psi-n+1-hn}
\Psi(n+1)+\gamma=H_n,
\end{equation}
in which $\gamma=0.57721\cdots$ is the Euler constant
\cite{handbook}. Therefore, we obtain
\begin{equation}\label{l-summing-psi}
\sum_{k=1}^n\left\{
\frac{(\Psi(k+1)+\gamma)^2}{k}-\frac{\Psi(k+1)+\gamma}{k^{2}}\right\}=\frac{(\Psi(n+1)+\gamma)^3-\zeta_n(3)}{3}.
\end{equation}
Letting
$$
\S(m,n)=\sum_{k=1}^n\frac{\Psi(k)^m}{k},
$$
the following identity \cite{hassani} is a result of 2-dimensional
$L-$summing method
\begin{equation}\label{l-summing-psik/k}
\S(1,n)=\frac{(\Psi(n+1)+\gamma)^2+\Psi(1,n+1)}{2}-\frac{\pi^2}{12}-\Psi(
n+1)\gamma-\gamma^2,
\end{equation}
where $\Psi(m,x)=\frac{d^m}{dx^m}\Psi(x)$ is called $m^{th}$
polygamma function \cite{handbook} and we have
$\sum_{k=1}^n\frac{1}{k^2}=-\Psi(1,n+1)+\frac{\pi^2}{6}$, which is
special case of the following identity
\begin{equation}\label{zeta-n-polygamma}
\zeta_n(s)=\frac{(-1)^{s-1}}{(s-1)!}\Psi(s-1,n+1)+\zeta(s)\hspace{10mm}(s\in\mathbb{Z},
s\geq 2),
\end{equation}
and using it in (\ref{l-summing-zeta-n}) one can get a
generalization of (\ref{l-summing-psi}), however this relation
itself is the key of getting an analogue of (\ref{l-summing-psik/k})
in $\R^3$, stated bellow.
\begin{thm} For every integer $n\geq 1$, we have
$$
\sum_{k=1}^n\frac{\Psi(k)^2}{k}+\frac{\Psi(k)}{k^2}=\frac{(\Psi(n+1)+\gamma)^3}{3}
-\frac{\zeta_n(3)}{3}+(\gamma-2)\frac{\pi^2}{6}-(\gamma-2)\Psi(1,n+1)-\gamma^2\Psi(n+1)-\gamma^3-2\S(1,n).
$$
\end{thm}
\begin{proof} We begin from the left hand side of the identity
(\ref{l-summing-psi}); simplifying it by the relations
(\ref{psi-n+1-n}), (\ref{psi-n+1-hn}) and the relation
(\ref{zeta-n-polygamma}) with $s=2$, gives the result.
\end{proof}

\begin{cor} For every integer $n\geq 1$, we have
$$
\S(2,n)=\frac{(\Psi(n+1)+\gamma)^3}{3}
-\frac{\zeta_n(3)}{3}+(\gamma-2)\frac{\pi^2}{6}-(\gamma-2)\Psi(1,n+1)
-\gamma^2\Psi(n+1)-\gamma^3-2\S(1,n)-\sum_{k=1}^n\frac{\Psi(k)}{k^2}.
$$
\end{cor}
In above corollary, the main term in the right hand side is
$\frac{\Psi(n+1)^3}{3}$. Also, we note that the summation
$\sum_{k=1}^n\frac{\Psi(k)}{k^2}$ is converges. Thus, we can write
the following asymptotic relation
$$
\S(2,n)=\frac{\Psi(n+1)^3}{3}+O(\ln^2n)\hspace{10mm}(n\rightarrow\infty).
$$
Similarly, considering (\ref{l-summing-psik/k}) we have
$$
\S(1,n)=\frac{\Psi(n+1)^2}{2}+O(\ln
n)\hspace{10mm}(n\rightarrow\infty).
$$
\textit{Note and Problem.} It is interesting to find an explicit
(probably recurrence) relation for the function $\S(m,n)$.
Considering two above asymptotic relations, we guess that
$$
\S(m,n)=\frac{\Psi(n+1)^{m+1}}{m+1}+O(\ln^m
n)\hspace{10mm}(n\rightarrow\infty).
$$
One can attack to this problem considering generalization of
$L-$summing method in higher dimension spaces, pointed in the last
section of this paper.

\section{Generating some new identities by Maple and 3-dimension $L-$summing
method} Appendix of this paper includes Maple program of 3-dimension
$L-$summing method. By $\textit{\texttt{LSMI}}<A_{abc}>$, we call
the identity outputted by $L-$Summing method's Maple program with
input $A_{abc}$. The algorithm of this program is result of the
formulation of 3-dimension $L-$summing method in above sections. In
this program we input a 3-dimensional array $A_{abc}$, then out put
is an identity generated by Maple. In this section we will state
some of these identities, with handling a detailed proof.
\subsection{Some elementary functions}
\begin{prop}\label{lsmi-ln-a} We have
$$
\texttt{LSMI}<\ln(a)>: \sum _{k=1}^{n}\left\{{k}^{2}\ln
k+2\,k\ln\Gamma(k+1)-2\,k\ln k-\ln\Gamma(k+1)+\ln k \right\}
={n}^{2} \ln \Gamma(n+1).
$$
\end{prop}
\begin{proof} Considering the relations (\ref{3d-l-summing}) and $\Gamma(n+1)=n!$,
we have $\Sigma(n)=n^2\sum_{a=1}^n\ln a=n^2\ln\Gamma(n+1)$. Also,
$\Sigma_2=k^2\ln k+2k\ln\Gamma(k+1)$, $\Sigma_1=\ln\Gamma(k+1)+2k\ln
k$ and $\Sigma_0=\ln k$. Putting these relations in
(\ref{3d-l-summing}) yields $\textit{\texttt{LSMI}}<\ln(a)>$ as
desired.
\end{proof}
\begin{cor} We have
$$
\sum_{k=1}^n \left\{(k^2-k)\ln
k+2k\ln\Gamma(k+1)\right\}=(n^2+n)\ln\Gamma(n+1).
$$
\end{cor}
\begin{proof} Breaking up the statement under the summation obtained from
$\textit{\texttt{LSMI}}<\ln(a)>$ in Proposition \ref{lsmi-ln-a},
into the sum of $(k^2-k)\ln k+2k\ln\Gamma(k+1)$ and $\ln\Gamma
\left( k+1 \right)+k\ln k -\ln k$, and considering the Proposition 6
of \cite{hassani}, which states
$$
\sum _{k=1}^{n}\left\{\ln\Gamma \left( k+1 \right)+k\ln k -\ln
k\right\}=n\ln \Gamma(n+1),
$$
completes the proof.
\end{proof}
\begin{rem} Examining Maple code of expressed summation on above
corollary, one can see that Maple has no comment on the computing
this summation; however, it is obtained by Maple itself and
$L-$summing method. This example shows that program-writers of Maple
can add $L-$summing method in the summation package of this
software, in order to making it able to compute some summations
which already couldn't compute them.
\end{rem}

\begin{prop} A little simplifying $\texttt{LSMI}<\tan(a)>$, we have
$$
\sum_{k=1}^n\left\{(k-1)^2\tan k+(2k-1)\T(k)\right\}=n^2\T(n),
$$
where $\T(n)=\sum_{k=1}^n\tan k$.
\end{prop}
\begin{proof} Considering the relation (\ref{3d-l-summing}),
we have $\Sigma(n)=n^2\sum_{a=1}^n\tan a=n^2\T(n)$. Also,
$\Sigma_2=k^2\tan k+2k\T(k)$, $\Sigma_1=\T(k)+2k\tan k$ and
$\Sigma_0=\tan k$, and consequently $L_k=(k-1)^2\tan k+(2k-1)\T(k)$.
This completes the proof.
\end{proof}
\subsection{Hypergeometric functions} In the next proposition, we introduce an identity
concerning hypergeometric functions, denoted in Maple by
$\texttt{hypergeom([}a_1~~a_2~\cdots~a_p\texttt{],[}b_1~~b_2~\cdots~b_q\texttt{],}x\texttt{)}$.
Standard notation and definition \cite{pet-wil-zeil} is as follows
$$
~_pF_q\left[\begin{array}{cccc}
  a_1 & a_2 & \cdots & a_p \\
  b_1 & b_2 & \cdots & b_q
\end{array};x\right]=\sum_{k\geq 0}t_k
x^k,
$$
where
$$
\frac{t_{k+1}}{t_k}=
\frac{(k+a_1)(k+a_2)\cdots(k+a_p)}{(k+b_1)(k+b_2)\cdots(k+b_q)(k+1)}x.
$$
\begin{prop} A little simplifying $\texttt{LSMI}<a!>$ and stating in standard notations, we have
$$
\sum_{k=1}^n\left\{
(k-1)^2k!+(2k-1)(k+1!)\H(1,k+2)\right\}=n^2(n+1)!\H(1,n+2),
$$
where
$$
\H(\alpha,\beta)=~_2F_0\left[\begin{array}{cc}
  \alpha & \beta \\
  -
\end{array};1\right].
$$
\end{prop}
\begin{proof} Considering definition of hypergeometric functions we
have $\H(1,n+1)=(n+1)\H(1,n+2)$, which implies
$\sum_{a=1}^na!=\H(1,2)-(n+1)!\H(1,n+2)=\P(n)$, say. This gives
$\Sigma(n)=n^2\P(n)$ and in similar way it yields that
$L_k=(k-1)^2k!+(2k-1)\left((k+1!)\H(1,k+2)-\H(1,2) \right)$. Thus,
we obtain
$$
\sum_{k=1}^n\left\{ (k-1)^2k!+(2k-1)\left((k+1!)\H(1,k+2)-\H(1,2)
\right)\right\}=n^2(n+1)!\H(1,n+2)-n^2\H(1,2),
$$
and a easy simplifying this, implies the result.
\end{proof}
\begin{rem} Three last propositions are examples of the array
$A_{abc}=f(a)$, for some given function $f$. In this case,
$L-$summing method takes the following form
$$
\sum_{k=1}^n\left\{(2k-1)\F(k)+(k-1)^2f(k)\right\}=n^2\F(n),
$$
where $\F(n)=\sum_{a=1}^n f(a)$.
\end{rem}

\section{Further generalizations of the $L-$summing method and some comments}

\subsection{The $L-$Summing method in $\R^t$} Consider a $t-$dimensional array $A_{x_1x_2\cdots x_t}$
and let $\Sigma(n)=\sum A_{x_1x_2\cdots x_t}$ with $1\leq
x_1,x_2,\cdots,x_t\leq n$. The $L-$Summing method in $\R^t$ is the
rearrangement $\Sigma(n)=\sum L_k$, where
$L_k=\sum_{m=1}^t\left\{(-1)^{m-1}\Sigma_{t-m}\right\}$ and
$$
\Sigma_{t-m}=\sum_{1\leq i_1<i_2<\cdots<i_m\leq
t}\left\{{\sum}^{\prime}A_{\textbf{x}_{i_1i_2\cdots i_m}}\right\},
$$
where in the inner summation ${\sum}^{\prime}$ is over
$x_j\in\{x_{i_1},\cdots,x_{i_m}\}^{C}=\{x_1,x_2,\cdots,x_t\}-\{x_{i_1},\cdots,x_{i_m}\}$
with $1\leq x_j\leq k$, and the index $\textbf{x}_{i_1i_2\cdots
i_m}$ denotes $x_1x_2\cdots x_t$ with
$x_{i_1}=x_{i_2}=\cdots=x_{i_m}=k$. One can apply this generalized
version to get more general form of relations obtained in previous
sections. For example, considering the array $A_{x_1x_2\cdots
x_t}=(x_1x_2\cdots x_t)^{-s}$ with $s\in\mathbb{C}$, yields
$$
\sum_{k=1}^n\left\{\sum_{m=1}^{t-1}(-1)^{m-1}{t\choose
m}k^{-ms}\zeta_{k}(s)^{t-m}\right\}=\zeta_n(s)^t+(-1)^t\zeta_n(ts).
$$

\subsection{$L-$summing method on manifolds} As in rising of this
paper, the base of the $L-$summing method is ordinary
multiplications table. Above generalization of the $L-$Summing
method in $\R^t$ is based on generalized multiplication tables
\cite{hassani-mt}. But, $\mathbb{R}^t$ is a very special
$t$-dimensional manifold, and if we replace it with $\Gamma$, an
$l-$dimensional manifold with $l\leq t$, then we can define
generalized multiplication table on $\Gamma$ by considering lattice
points on it (which of course isn't easy problem). Let
$$
L_\Gamma(n)=\big\{ (a_1, a_2,\cdots, a_t)\in\Gamma\cap\mathbb{N}^t:
1\leq a_1,a_2,\cdots,a_t\leq n \big\},
$$
and $f:\mathbb{R}^k\longrightarrow\mathbb{C}$ is a function. If
$\mathcal{O}_{\Gamma}$ is a collection of $k-1$ dimension orthogonal
manifolds, in which
$L_\Gamma(n)=\cup_{\Lambda\in\mathcal{O}_{\Gamma}}L_{\Lambda}(n)$
and $L_{\Lambda_i}(n)\cap L_{\Lambda_j}(n)=\phi$ for distinct
$\Lambda_i,\Lambda_j\in\mathcal{O}_{\Gamma}$, then we can formulate
$L-$summing method as follows,
$$
\sum_{X\in
L_\Gamma(n)}f(X)=\sum_{\Lambda\in\mathcal{O}_{\Gamma}}\left\{\sum_{X\in
L_{\Lambda}(n)}f(X)\right\}.
$$
Here $L-$summing elements are $\sum_{X\in L_{\Lambda}(n)}f(X)$. This
may be useful when one apply it on some special manifolds.

\subsection{Stronger form of $L-$summing method} One can state the method of
$L-$summing $\sum L_k=\Sigma(n)$ in the following stronger form
$$
L_n=\Sigma(n)-\Sigma(n-1).
$$
Specially, this will be useful for those arrays with $\Sigma(n)$
computable explicitly and $L_k$ maybe note. For example, applying
this note on the array $A_{x_1x_2\cdots x_t}=(x_1x_2\cdots
x_t)^{-s}$ in $\R^t$ with $s\in\mathbb{C}$, implies
$$
\sum_{m=1}^{t-1}(-1)^{m-1}{t\choose
m}n^{-ms}\zeta_{n}(s)^{t-m}=\zeta_n(s)^t+(-1)^t\zeta_n(ts)-\zeta_{n-1}(s)^t-(-1)^t\zeta_{n-1}(ts).
$$

\bigskip
\hrule
\bigskip

\begin{center}
\textsc{Appendix. Maple Program of 3-dimension $L-$Summing Method
for the array $A_{abc}=\frac{1}{abc}$}
\end{center}
\small{\texttt{restart:\\
A[abc]:=1/(a*b*c);\\
S21:=sum(sum(eval(A[abc],a=k),b=1..k),c=1..k):\\
S22:=sum(sum(eval(A[abc],b=k),a=1..k),c=1..k):\\
S23:=sum(sum(eval(A[abc],c=k),a=1..k),b=1..k):\\
S2:=S21+S22+S23:\\
S11:=sum(eval(eval(A[abc],a=k),b=k),c=1..k):\\
S12:=sum(eval(eval(A[abc],a=k),c=k),b=1..k):\\
S13:=sum(eval(eval(A[abc],b=k),c=k),a=1..k):\\
S1:=S11+S12+S13:\\
S0:=eval(eval(eval(A[abc],a=k),b=k),c=k):\\
L[k]:=simplify(S2-S1+S0):\\
ST(A):=(simplify(sum(sum(sum(A[abc],a=1..n),b=1..n),c=1..n))):\\
Sum(L[k],k=1..n)=ST(A);}}
$$
A_{abc}:=\frac{1}{abc}
$$
$$
\sum _{k=1}^{n}{\frac {3\, \left( \Psi \left( k+1 \right)  \right)
^{2 }{k}^{2}+6\,\Psi \left( k+1 \right)
{k}^{2}\gamma+3\,{\gamma}^{2}{k}^{ 2}-3\,\Psi \left( k+1 \right)
k-3\,\gamma\,k+1}{{k}^{3}}}= \left( \Psi
 \left( n+1 \right) +\gamma \right) ^{3}
$$

\bigskip
\hrule
\bigskip

\end{document}